\documentclass[a4paper]{article}
\usepackage{graphicx, xcolor} 
\usepackage{amsmath, amsthm, amsfonts, amssymb, amscd, bm, enumerate}
\usepackage{verbatim}

\usepackage{color}
\usepackage{tikz}

\usepackage{url}

\renewenvironment{proof}{{\bfseries Proof:}}{\qed}
\newtheorem{theorem}{Theorem}

\newtheorem{lemma}{Lemma}
\newtheorem{definition}{Definition}
\newtheorem{assumption}{Assumption}
\newtheorem{remark}{Remark}

\def\F{\mathcal{F}}
\def\E{\mathcal{E}}

\def\bbH{\bar{H}}

\def\L{\mathcal{L}}
\def\bL{\bar{L}}

\def\U{\mathcal{U}}

\def\bR{\mathbb{R}}
\def\bE{\mathbb{E}}
\def\bN{\mathbb{N}}
\def\bu{\mathbf{u}}
\def\bbu{\bar{u}}
\def\bbv{\bar{v}}

\def\y{\tilde{y}}
\def\bz{\mathbf{z}}
\def\bl{\mathbf{\lambda}}
\def\bZ{\mathbf{Z}}
\def\Y{\tilde{Y}}
\def\Z{\tilde{Z}}

\def\bH{\mathbf{H}}
\def\bone{\mathbf{1}}
\def\bP{\mathbb{P}}

\def\bQ{\mathbb{Q}}

\mathchardef\mhyphen="2D

\def\<{<}
\def\>{>}

\begin{document}
\title{Nash Equilibria for non Zero-Sum Ergodic Stochastic Differential Games}
\author{Samuel N. Cohen\footnote{Research supported by the Oxford--Man Institute for Quantitative Finance and the Oxford--Nie Financial Big Data Laboratory.}\\ University of Oxford \\ \\ Victor Fedyashov\\ University of Oxford}

\date{\today}

\maketitle

\begin{abstract}
\noindent We consider non zero-sum games where multiple players control the drift of a process, and their payoffs depend on its ergodic behaviour. We establish their connection with systems of Ergodic BSDEs, and prove the existence of a Nash equilibrium under the generalised Isaac's conditions. We also study the case of interacting players of different type.

\noindent MSC: 91A60, 49N70, 93E20

\noindent Keywords: Ergodic games, Nash equilibrium, BSDE

\end{abstract}

\section{Introduction}
The area of non-zero sum stochastic differential games has been a subject of intensive research over the past several decades. It deals with the situation where multiple players are each trying to maximise their payoff, but rewards of all players do not sum up to any constant. The basic aim is to find a Nash equilibrium point -- a set of strategies for all players, such that none of them would benefit from unilaterally deviating. The study of the connection between this problem and backward stochastic differential equations (BSDEs) was pioneered by Hamad\`ene, Lepeltier and Peng in \cite{Karoui_BSDE}, building up on their previous work on zero-sum games in Hamad\`ene and Lepeltier \cite{HL}. Since then it has been discussed in multiple contexts, e.g. games of control and stopping by Karatzas and Li in \cite{KL} and risk-sensitive control by El Karoui and Hamad\`ene in \cite{KH}. The main goal of the present paper is to prove the existence of an equilibrium in the case where we have two drift controlling players who are ergodic, i.e. they choose their strategies in infinite horizon and optimise the long run average. In a standard fashion (see, e.g. Hamad\`ene and Mu \cite{Ham_Mu}), we impose the so-called generalised Isaac's condition on the Hamiltonian, in order to ensure the attainability of the infimum simultaneously for both players. We then prove that the existence of a saddle point for the game follows from the existence of a solution to a system of Ergodic BSDEs with continuous coefficients. Using a modified version of Picard iteration (not dissimilar in nature to the fixed point construction used by Hu and Tang in \cite{Hu_games}), along with an array of estimates for solutions of EBSDEs, established in Debussche, Hu and Tessitore \cite{Hu} and by the authors in \cite{Levy_paper}, we prove that such a system does admit a solution. We then extend these results to the case of finite or countably many players. As a byproduct of developing the machinery required to deal with ergodic games, we obtain the existence of a solution to a class of EBSDEs with continuous drivers of linear growth. The concluding section is dedicated to asymmetric games, where one player is purely ergodic, and the other also optimises in infinite horizon, but values the present more than the future (i.e. uses discounting). We show that the value of such a game exists, and demonstrate convergence to the ergodic system as we send the discount rate to zero.

\section{Setup}

We begin by describing the uncontrolled forward process $\{X_t\}_{t \geq 0}$, the drift of which will be affected by the players in the sequel. Let the dynamics of $\{X_t\}_{t \geq 0}$ be governed by the It\^o SDE:
\begin{equation}
	dX_t = (AX_t + F(X_t))dt + \sigma dW_t, \quad X_0 = x_0\in\mathbb{R}^N,
\label{forward_SDE}
\end{equation}
where $W$ denotes a standard Brownian motion in $\mathbb{R}^N$ under the measure $\bP$. We need a few preliminary assumptions.
\begin{assumption} Conditions to guarantee existence of a unique strong solution to (\ref{forward_SDE}):
\begin{itemize}
\item The triple $(\Omega, \F, \bP)$ is a probability space, $\{\F_t\}_{t\geq 0}$ is a right-continuous filtration and $\F_t$ is complete for each $t$.
\item The process $\{ W_t \}_{t \geq 0}$ is a Brownian motion with the predictable representation property in the filtration $\{\F_t\}_{t \geq 0}$.
\item The operator $A$ is dissipative and generates a stable semigroup $\{e^{tA}\}_{t \geq 0}$. In other words, there exist constants $\mu > 0$, $M > 0$ such that
\[
	\langle Ax,x \rangle \leq -\mu \| x\|^2 \text{  and  } \| e^{ t A} \| \leq M e^{-\mu t}
\]
hold for all $x \in \bR^n$ and all $t \geq 0$. Here and for the rest of the paper $\|\cdot \|$ denotes the Euclidean norm (or Frobenius norm for a matrix).
\item The map $F$ is uniformly Lipschitz continuous and bounded. 
\item The matrix $\sigma$ is invertible. 
\end{itemize}
\begin{remark}
 As we will be considering the ergodic behaviour of (\ref{forward_SDE}), under a variety of probability measures, we need to be careful not to complete our filtrations in the usual way. The above requirement is non-standard, but it is easy to see that all the usual existence proofs hold via pasting in time.
\end{remark}

\label{existence_solution}
\end{assumption}
 The following result is a straightforward application of Gr\"onwall's lemma and the dissipativity condition:
\begin{theorem} If Assumption \ref{existence_solution} is satisfied, then there exists a unique strong solution to the equation (\ref{forward_SDE}). Moreover,
\[
	\bE \big(  \| X_t \|^2 \big) < C (1 + \|x_0\|^2),
\]
where the constant $C$ depends on $M,\mu$, but is independent of $t$.
\label{E1}
\end{theorem}
\begin{proof} Given the assumption of Lipschitz continuity, the existence of a strong, square integrable solution to \eqref{forward_SDE} is standard. Fix a constant $k > 0$. Applying It\^o's formula to $e^{kt}\| X_t\|^2$ and taking expectations, we obtain
\[
	e^{kt} \bE \| X_t \|^2 = \bE \|x_0\|^2 + \int_0^t ke^{ks}\|X_s\|^2ds + \bE \int_0^t e^{ks}G(s)ds,
\]
where
\[
\begin{split}
	G(s) &= 2\langle AX_t + F(X_t), X_s  \rangle + \| \sigma \|^2 \\
		& \leq -2\mu \| X_s\|^2 + \langle F(X_s), X_s \rangle + \| \sigma\|^2.
\end{split}
\]
Taking into account Assumption \ref{existence_solution} and applying Cauchy--Schwarz, we see that for any $\epsilon > 0$, there exists $C^{\epsilon}$, such that 
\[
	G(s) \leq -2(\mu - \epsilon)\| X_s \|^2 + C^{\epsilon},
\]
for $0 \leq s \leq t$. Then, choosing $\epsilon < \mu$ and $k < 2(\mu - \epsilon)$, we arrive at
\[
	e^{kt}\bE \| X_t \|^2 \leq \bE \| x_0 \|^2 + \frac{C^{\epsilon}}{k}e^{kt},
\]
concluding the proof. 

\end{proof}

\begin{remark}\label{remark2}By applying It\^o's lemma to $\|X_t\|^4$ and using the same reasoning as in the proof of Theorem \ref{E1}, one can also show that 
\begin{equation}
	 \bE \| X_t \|^4 \leq C(1 + \| x_0 \|^4),
\end{equation}
for some constant $C$ independent of time. 
\end{remark}

\begin{remark} If $F(\cdot)$, the nonlinear part of the drift, is a bounded measurable map, there still exists a solution to (\ref{forward_SDE}) in the weak sense. In other words, for any $T>0$, there exists\footnote{
If we assume that $\Omega$ is a canonical Wiener space, then Kolmogorov's extension theorem can be used to define a measure $\tilde \bP$ such that $\tilde \bP|_{\F_T} = \tilde\bP^T|_{\F_T}$ for all $T>0$, but $\tilde\bP$ is generally not equivalent to $\bP$. This is why we did not assume the usual conditions on $\{\F_t\}_{t\ge 0}$, as the presence of too many null sets in $\F_0$ would be problematic. In more general settings, this extension is a nontrivial exercise, see, for example, Parthasarathy \cite{Parthasarathy}.
}  a measure $\tilde{\bP}^T\sim \bP$, and an $\tilde{\bP}^T$-Brownian motion $\tilde{W}$ such that, for $t\in[0, T]$,
\[
	X_t = e^{tA}x_0 + \int_0^t e^{sA}F(X_s)ds + \int_0^t e^{sA}\sigma d\tilde{W}(s),
\]
 and this solution is unique in law (and is consistent for different values of $T<\infty$, so $\tilde W$ is uniquely defined). Moreover, since this solution can be obtained through a measure change from a strong solution to the equation without drift, $\{ \tilde{W}_t \}_{t\in [0,T]}$ has the predictable representation property in the filtration $\{\F_t\}_{t \in [0,T]}$ under $\tilde \bP^T$, for all $T$. For details see, e.g. Chapters 15 and 18 in \cite{SCA}.
\end{remark}

 Let $\U_1 \times \U_2$ be a separable metric space of controls in which $(u_t(\omega), v_t(\omega))$ takes values. Denote by $L_i: \bR^N \times \U_1 \times \U_2 \to \bR$, $i=1,2$, bounded measurable cost functions, such that 
\[
	| L_i(x,u,v) - L_i(x',u,v)| \leq C \| x - x' \|,
\]
for some $C > 0$. Given an arbitrary pair of admissible controls $(u,v)$  (i.e. predictable stochastic processes $(\{u_t\}_{t \geq 0}, \{v_t\}_{t \geq 0})$), we define the Girsanov density 
\[
	\rho^{u,v}_T := \exp \bigg( \int_0^T R(u_t,v_t)dW_t - \frac{1}{2}\int_0^T \| R(u_t,v_t) \|^2 ds \bigg)
\]
for some function $R:\U_1\times \U_2 \to \bR^N$, satisfying $\| R(\cdot,\cdot) \| \leq \bar{C}$ for some $\bar{C} > 0$. We define the probability $\bP^{u,v,T} := \rho^{u,v}_T\bP$. For each agent ($i = 1,2$) we define ergodic payoffs 
\begin{equation}
	J^i(u, v) = {\lim\sup}_{T\to\infty} T^{-1} \bE^{u,v,T}\bigg[\int_0^T L_i(X_t, u_t,v_t) dt\bigg].
\label{payoff}
\end{equation}
The goal is to find an admissible control $(u^*,v^*)$, such that 
\[
	J^1(u^*,v^*) \leq J^1(u,v^*) \text{  and  } J^2(u^*,v^*) \leq J^2(u^*,v)
\]
holds for all admissible $(u,v)$. 

 In order to use dynamic principles to determine optimal $u$ and $v$, we define the Hamiltonian functions $\{H_i\}_ {i = 1,2}$, as 
\begin{equation}
	H_i(x,z_i,u,v) = z_i R(u,v) + L_i(x,u,v)
\label{H}
\end{equation}
for $(u,v) \in \U_1 \times \U_2$. We notice that for a fixed pair of controls $(u,v)$, for each player ($i=1,2$), the payoff (\ref{payoff}) is an ergodic average. Consider the following ergodic backward stochastic differential equation (EBSDE), in which the Hamiltonian plays the role of the driver term:
\begin{equation}
	Y^{i}_t = Y^i_T+\int_t^T[H_i(X_{s}, Z^i_s,u_s,v_s)-\lambda^{i}] ds - \int_t^T Z^i_s dW_s.
\label{Etmp}
\end{equation}
Suppose there exists a solution triplet $(Y^i,Z^i,\lambda^i)$ to (\ref{Etmp}), and there also exists a constant $C = C(x_0) > 0$ independent of $T$, such that $\bE^{u,v,T}|Y_t^i| < C$ for all $0 \leq t \leq T<\infty$. Then, by changing measure and dividing by $T$, we immediately see that 
\[
	J^i(u, v) = \lambda^i.
\]
We refer to (\ref{Etmp}) as the BSDE corresponding to the value function (\ref{payoff}). In order to proceed we need the following assumptions:

\begin{assumption} Part (i) below is called the generalised Isaac's condition. Part (ii) is a technical condition, the motivation behind it will be discussed later. 
\begin{enumerate}[(i)]
 \item There exist measurable maps $u_1^*,u_2^*$, defined on $ \bR^{3N}$, with values in $\U_1$ and $\U_2$ respectively, such that 
\[
	H_1(x,z_1,u_1^*(x,z_1,z_2),u_2^*(x,z_1,z_2)) \leq H_1(x,z_1,u,u_2^*(x,z_1,z_2))
\]
 and 
 \[
	H_2(x,z_1,u_1^*(x,z_1,z_2),u_2^*(x,z_1,z_2)) \leq H_2(x,z_1,u_1^*(x,z_1,z_2),v) 	
 \]
 holds for all $(x,z_1,z_2,u,v) \in  \bR^{3N} \times \U_1 \times \U_2$. 
 \item The mapping $(z_1,z_2) \in \bR^{2N} \to (H^*_1,H^*_2)(x,z_1,z_2)$ is continuous for any fixed $x$. 
\end{enumerate}
\label{Isaacs}
\end{assumption}

 One way to think about the condition above is as the natural extension of the standard assumptions used in the zero-sum framework to the non zero-sum case. Indeed, consider the case of two players, one maximising and the other minimising a functional 
\[
	J(u,v) = {\lim\sup}_{T\to\infty} T^{-1} \bE^{u,v,T}\bigg[\int_0^T L(X_t, u_t,v_t) dt\bigg],
\]
by choosing controls $(u_t,v_t) \in \U_1 \times \U_2$ for $t \geq 0$. Define the Hamiltonian as $H(x,z,u,v) = zR(u,v) + L(x,u,v)$. Then, using a measurable selection argument, one could show that the standard Isaac's condition
\[
	\max_{v \in \U_2} \min_{u \in \U_1} H(x,z,u,v) = \min_{u \in \U_1} \max_{v \in \U_2} H(x,z,u,v), \quad \forall x,z \in \bR^N
\]
is equivalent to the existence of Borel-measurable functions $u^*:\bR^{2N} \to \U_1$ and $v^*:\bR^{2N} \to \U_2$, such that 
\[
\begin{split}
	&H(x,z,u^*(x,z),v^*(x,z)) \leq H(x,z,u,v^*(x,z)) \text{  for all  } u \in \U_1, \\
	&H(x,z,u^*(x,z),v^*(x,z)) \geq H(x,z,u^*(x,z),v) \text{  for all  } v \in \U_2. \\
\end{split}
\]
For details, see Hamad\`ene and Lepeltier \cite{HL}. Given this connection, one can see that (i) in the Assumption \ref{Isaacs} is a reasonable generalisation of the Isaac's condition. 

\begin{remark} We notice that, for fixed $u,v$, the function $H_i$ is Lipschitz in $z_i$. However, Assumption \ref{Isaacs} only gives us the existence of a so called ``closed loop'' control, where the optimal processes $u^*,v^*$ depend on $(z_1,z_2)$, and only continuity of $H_i$ is guaranteed. Therefore the standard method of proving the existence of a solution to the corresponding Ergodic BSDE breaks down. %ZERO SUM%

We should observe that Assumption \ref{Isaacs} does not guarantee that either $u^*$ or $v^*$ are continuous (cf. Mannucci \cite{Mannucci} for similar assumptions).
\end{remark}

\begin{remark} In the game we consider, both players control the drift. In \cite{Ham_Mu} the authors considered the case where the resulting $R(u,v)$ can be unbounded, but is of linear growth with respect to the underlying process $\{X_t\}_{t \geq 0}$. However, since we are in the ergodic framework, we restrict our attention to the bounded case.  While Isaac's condition is arguably natural, we should note that there has been work on stochastic differential games without this assumption (see for example Buckdahn, Li and Quincampoix \cite{BuckdahnLi}).
\end{remark}

\begin{remark} One may also notice the similarity between the generalised Isaac's condition and the definition of a saddle point. Indeed, one could think of part (i) of Assumption \ref{Isaacs} as ensuring the existence of Nash equilibrium controls on the infinitesimal scale. 

\label{remarkTmp}
\end{remark}

%\subsection{Zero sum games}

\section{Ergodic BSDEs with continuous coefficients and stochastic games}

In this section we prove that, under our assumptions, the game has value, that is there exists a Nash equilibrium as defined in the preceding section. We begin by considering two players of ergodic type, then proceed to generalise to the case of $n \in \bN \cup \{+\infty\}$ players. Finally, we show a connection with PDEs. 

\subsection{Games with two players of ergodic type}

 In order to establish the existence of a Nash equilibrium, the following two-step programme needs to be carried out: 
\begin{itemize} 
\item Assuming there exists a solution to the EBSDE corresponding to control pair $(u^*, v^*)$ as defined in Assumption \ref{Isaacs}, establish its optimality. 
\item Prove the existence of a unique solution to the EBSDE with continuous coefficients, guaranteeing this assumption holds.
\end{itemize}

 In order to establish the existence of a solution to a one-dimensional ergodic BSDE with driver $f$, the following set of assumptions on the driver needs to be made:

\begin{assumption} The `driver' $f:\bR^N \times \bR^N \to \bR$ is a measurable map. Moreover, there exist constants $l,\kappa$ such that 
\[
	\big| f(x,0) \big| \leq l, \quad \big| f(x,z) - f(x,\bar{z}) \big| \leq \kappa \| z - \bar{z} \|
\]
holds for all $x,z,\bar{z} \in \bR^N$.
\label{A1}
\end{assumption}

 We begin by stating certain crucial results from the theory of Ergodic BSDEs. The following theorem establishes existence and uniqueness of solutions. We omit the proof, as it can be easily constructed using the arguments in \cite{Hu}.

\begin{theorem} Let the driver $f$ satisfy Assumption \ref{A1}. Then there exists a solution $(Y,Z,\lambda)$, where $\{Y_t\}_{t\geq 0}$ is adapted, $\{Z_t\}_{t \geq 0}$ is predictable, and $\lambda \in \bR$ is a constant, to the EBSDE
\[
	 Y_t = Y_T+\int_t^T[f(X_s,Z_s)-\lambda] ds - \int_t^T Z_s dW_s,	
\]
such that $| Y_t | \leq C (1 + \| X_t \|^2)$ for some constant $C>0$. There also exists a deterministic function $v:\bR^N \to \bR$, such that $Y_t = v(X_t)$. Moreover, the solution is unique in the class of Markovian solutions, for which the $Y$ component is of polynomial growth with respect to $X$. 
\label{TErg}
\end{theorem}

\if 0

\begin{remark}
 We should note that the existence results in \cite{Hu} assume that the volatility in
the forward dynamics is constant with respect to $X_t$. The critical point where this is needed in \cite{Hu} is in the proof of the ergodicity of the forward process for the class of drift with bounded nonlinear part. Since in the present paper we are dealing with a finite-dimensional process, this proof (and thus existence of a solution to the EBSDE) can be adapted to incorporate state-dependence if the following way. We first show (for proof see Theorem 7 in \cite{Adjoint_paper}) that if the nonlinear part of the drift $F$ is Lipschitz, then there exist constants $C > 0$ and $\rho > 0$ such that, for any bounded continuous function $\psi : \bR^n \to \bR$,
\[
	|P(\tau,t)[\psi](x)  - P(\tau,t)[\psi](y) | \leq C (1 + ||x||^2 + ||y||^2)e^{-\rho (t-\tau)}\|\psi(u)\|_0,
\]
where $P_t \psi (x) := \bE \big[ \psi (X^x_t) \big]$. By Corollary 2.5 in \cite{Hu} we can extend this result to the class of bounded drifts that can be approximated by a uniformly bounded sequence of Lipschitz functions. What remains is to notice that when trying to prove the existence of a solution to an ergodic BSDE, we only deal with additional drifts that come from change of measure, and are thus, by Lemma 3.5 in \cite{Hu}, of the required class.

\end{remark}

\fi

 We now show that if both players can find a Markovian solution of polynomial growth to their respective EBSDEs with optimal pair of controls $(u^*,v^*)$, then by so doing they find a Nash equilibrium. 

\begin{theorem} Assume that there exist two triplets $(Y^i, Z^i, \lambda^i)$, $i=1,2$, such that 
\begin{equation}
	 Y^{i}_t = Y^i_T+\int_t^T[H_i(X_{s}, Z^i_s,u^*(Z^1_s,Z^2_s),v^*(Z^1_s,Z^2_s))-\lambda^{i}] ds - \int_t^T Z^i_s dW_s,	
\label{EBSDE}
\end{equation}
holds for all $0 \leq t \leq T < 0$, where the pair $(u^*,v^*)$ is as in Assumption \ref{Isaacs}. Moreover, assume there exist two deterministic functions with polynomial growth $y_i(x)$, $i=1,2$, such that $Y^i_t = y_i(X_t)$ holds $\bP$-a.s. for all $t \geq 0$. Then the control $(u^*(X_s,Z^1_s,Z^2_s),v^*(X_s,Z^1_s,Z^2_s))$ is a Nash equilibrium.
\label{T3}
\end{theorem}
 \textbf{Proof:} We focus on the case $i=1$. For $t \geq 0$ write $u^*_t = u^*(Z^1_t,Z^2_t)$ and $v^*_t = v^*(Z^1_t,Z^2_t)$. The control $(u^*,v^*)$ is clearly admissible. For an arbitrary admissible control $u$, for any $T>0$, the pair $(u,v^*)$ generates a measure $\bP^{u,v^*,T}$, under which $W^{u,v^*}_t = W_t - \int_0^t R(u_s,v^*_s)ds$ is a Brownian motion on $[0,T]$. Then, since $(Y^1, Z^1, \lambda^1)$ solves the EBSDE,  
\[
\begin{split}
	\lambda^1 T = Y^1_T-Y^1_0 &+ \int_0^T [H_1(X_{s}, Z^1_s,u^*_s,v^*_s) - Z^1_sR(u_t,v^*_s) - L_1(X_s,u_s,v^*_s)]ds \\
	&- \int_0^T Z^1_s dW^{u,v^*}_s + \int_0^T L_1(X_s,u_s,v^*_s)ds.
\end{split}
\]
Taking expectations with respect to $\bP^{u,v^*,T}$, and remembering that $\lambda^1$ is a constant, we arrive at
\[
	\lambda^1 \leq \frac{1}{T} \bE^{u,v^*,T} [Y^1_T - Y^1_0] + \frac{1}{T}\bE^{u,v^*,T} \bigg[ \int_0^T L_1(X_s,u_s,v^*_s)ds \bigg],
\]
since 
\[
	H_1(X_{s}, Z^1_s,u^*_s,v^*_s) - Z^1_sR(u_t,v^*_s) - L_1(X_s,u_s,v^*_s) \leq 0
\]
by the definition of $(u^*,v^*)$. We recall that $Y^1$ is of polynomial growth in $X$, and the latter has all moments under $\bP^{u,v^*}$ uniformly in $t$ (for details see, e.g. \cite{Hu}). Thus, taking $\limsup$ on both sides, we obtain 
\[
	\lambda^1 \leq \limsup_{T \to \infty} \frac{1}{T}\bE^{u,v^*,T} \bigg[ \int_0^T L_1(X_s,u_s,v^*_s)ds \bigg] = J^1(u,v^*).
\]
For $u = u^*$ the inequality everywhere becomes equality by the definition of the pair $(u^*,v^*)$ in Assumption \ref{Isaacs}.  In the same way, we observe $J^2(u^*, v^*)=\lambda^2\leq J^2(u^*, v)$ for any admissible control $v$, finishing the proof.

\qed

\begin{remark} The result of Theorem \ref{T3} shows that even in the case where admissible controls are allowed to be non-Markovian, the values $\lambda^{i}$, $i=1,2$ obtained as solutions to the system (\ref{EBSDE}) are still optimal. In other words, using the generalised Isaac's condition, we can show the existence of ``closed loop'' or ``feedback'' controls that are optimal in the larger class of ``open loop'' or ``predictable'' controls. 
\end{remark}

 We now proceed to the second step of our programme. We prove that there exists a solution to the system of ergodic BSDEs (\ref{EBSDE}). The following auxiliary result is a weak version of a comparison theorem for the ergodic values $\lambda$, which holds in multiple dimensions. 

\begin{lemma} Let $f^i$, $i=1,2$ be two drivers, such that for any $(x,z_1,z_2)$
\[
	| f^i(x,z_1,z_2)| \leq C\|z_i\| + \bar{C}
\]
for some constants $C,\bar{C}>0$. Suppose that there exist Markovian solutions $(Y^i,Z^i,\lambda^i)$, $i=1,2$ to the corresponding EBSDEs. Let $\lambda$ be the ergodic part of the solution $(Y,Z,\lambda)$ to the EBSDE
\[
	 -dY_t = [(C\| Z_t \| + \bar{C})-\lambda] ds - Z_t dW_t.	
\]
Then $\lambda^i \leq \lambda$ for $i=1,2$.
\label{L1}
\end{lemma}

 \textbf{Proof: } Define a measure $\bQ^T$ as follows:
\[
	\frac{d\bQ^T}{d\bP} := \E \bigg( \int_0^T \frac{f(Z_t) - f(Z^i_t)}{\| Z_t-Z^i_t \|^2} ( Z_t-Z^i_t)dW_t \bigg),
\]
where $f(z) = C\| z \| + \bar{C}$. Then under $\bQ^T$, the process $\{X_t\}_{t \geq 0}$ is ergodic with some invariant measure $\mu$. Since we only consider Markovian solutions, we know that $Y^i_t = v^i(X_t)$, $Y_t = v(X_t)$ and $Z^i_t = \zeta^i(X_t)$ (for details on Markov representation, see e.g. Fuhrman, Hu and Tessitore \cite{Hu_Banach}). Then, integrating the difference $v^i(x)-v(x)$, we obtain
\[	
\begin{split}
	&\int_{\bR^N}(v^i(x)-v(x))\mu(dx) \\
	& = \int_{\bR^N}  \bE^{Q^T} \bigg[ (v^i(X_T)-v(X_T)) +  \int_0^T f^i (X_t,Z^1_t,Z^2_t)-f(Z^i_t) dt - T(\lambda^i - \lambda) \bigg] \mu(dx)  \\
	& = \int_{\bR^N}(v^i(x)-v(x))\mu(dx) - T(\lambda^i - \lambda)\\
 &\qquad  + \int_0^T \int_{\bR^N} f^i (x,\zeta^1(x),\zeta^2(x)) -f(\zeta^i(x))\mu(dx) dt,
\end{split}
\]
from which we immediately see that 
\[
	\lambda^i - \lambda = \int_{\bR^N} (f^i(x,\zeta^1(x),\zeta^2(x)) - f(\zeta^i(x))) \mu(dx) \leq 0,
\]
finishing the proof. 

\qed

 There are two apparent ways to attack the problem of proving the existence of a solution to (\ref{EBSDE}):
\begin{enumerate}[(i)] \item For any $t \geq 0$, equation (\ref{EBSDE}) gives us a map $Y^i_{t+1} \to Y^i_t$ for $i \in 1,2$. If we integrate with respect to the ergodic measure $\pi$ for the forward process $\{X_t\}_{t \geq 0}$, we have a new map. Now we know that the solution to (\ref{EBSDE}) is its fixed point. To proceed, we would need to prove certain compactness and continuity properties for this map.

\item We could attempt construction of a solution iteratively using a form of Picard iteration. 
\end{enumerate}

 In this paper we focus on plan (ii). The main result of this section is as follows:

\begin{theorem} Suppose Assumptions \ref{existence_solution}, \ref{Isaacs} and \ref{A1} are satisfied. Then there exist solution triplets $(Y^i, Z^i, \lambda^i)$, $i=1,2$, where $Y^i$ is adapted, $Z^i$ predictable, and $\lambda^i$ a constant, such that 
\[
	 Y^{i}_t = Y^i_T+\int_t^T[H_i(X_{s}, Z^i_s,u^*(Z^1_s,Z^2_s)(x),v^*(Z^1_s,Z^2_s)(x))-\lambda^{i}] ds - \int_t^T Z^i_s dW_s,	
\]
holds for all $0 \leq t \leq T < 0$.  Moreover, there exist deterministic functions with polynomial growth $y^1(x),y^2(x)$, such that $Y^i_t = y^i(X_t)$, $i=1,2$ holds $\bP$ - a.s. for all $t \geq 0$. 
\label{T4}
\end{theorem}

 \begin{proof} We begin by setting $Z^{i,0} \equiv 0$ for $i=1,2$. Given a set of  Markovian solutions $\{(Y^{i,n-1},Z^{i,n-1},\lambda^{i,n-1})\}_{i=1,2}$ obtained previously, there exist deterministic functions $\{ \zeta^{i,n-1} \}_{i=1,2}$, such that $Z^{i,n-1}_t = \zeta^{i,n-1}(X_t)$. We define the next iteration as a solution to the EBSDE
\begin{equation}
	 dY^{i,n}_t = - [f^{i,n}(X_t,Z^{i,n}_t)-\lambda^{i,n}]dt +  Z^{i,n}_t dW_t,		\quad i=1,2
\label{E2}
\end{equation}
where 
\[
	f^{i,n}(x,z) = H_i\Big(x, z,u^*\big(\zeta^{1,n-1}(x),\zeta^{2,n-1}(x)\big),v^*\big(\zeta^{1,n-1}(x),\zeta^{2,n-1}(x)\big)\Big). 
\]
We note that for all $n \geq 1$, the drivers $\{f^{i,n} \}_{i=1,2}$ are uniformly Lipschitz in the $z$ component, since 
\[
	\big| f^{i,n}(x,z) - f^{i,n}(x,z') \big| = \big\| R((u^*,v^*)(\zeta^{1,n-1},\zeta^{2,n-1})) \big\| \| z-z' \| \leq C \| z-z' \|.
\]
We therefore know that there exists a solution $\{(Y^{i,n},Z^{i,n},\lambda^{i,n}) \}_{i=1,2}$ to (\ref{E2}). 

It remains to show that there exists a convergent subsequence with limit $\{(Y^i,Z^i,\lambda^i)\}_{i=1,2}$. By Lemma \ref{L1}, there exists a constant $\lambda > 0$, such that 
\[
	| \lambda^{i,n} | \leq \lambda, \quad i = 1,2
\]
holds for all $n \geq 1$. We can therefore find a subsequence $\{n_k\}_{k \in \bN}$, such that $\lambda^{i,n_k} \to \lambda^i$, $i=1,2$ for some $\{ \lambda^i \}_{i=1,2}$. One can also show (see Theorem 8 in \cite{Levy_paper}), that for all $n \geq 1$, there exist deterministic functions $\{v^{i,n}\}_{i=1,2}$, such that $Y^{i,n}_t = v^{i,n}(X_t)$, $Z^{i,n}_t = \nabla v^{i,n}(X_t)\sigma$ and 
\begin{equation}
	| Y^{i,n}_t | = | v^{i,n}(X_t) | \leq C(1 + \| X_t \|^2)
\label{est}
\end{equation}
holds for some constant $C>0$ independent of $n$. (\textit{As shown in \cite{Levy_paper}, this constant depends on $\sup_{x \in \bR^n} F(x)$ and on the Lipschitz constant of $f^{i,n}(x,z)$ in the $z$ component. The fact that the latter is uniform in $n$ follows from the definition of the Hamiltonian.}) Moreover, we have the following gradient estimate:
\begin{equation}
	\| \nabla v^{i,n}(x) \| \leq C(1 + \| x \|^2) 
\label{grad}
\end{equation}
for $i=1,2,n \in \bN$. Since $\bR^N$ is separable, there exists a dense countable subset $B \subset \bR^n$. We can therefore extract a further subsequence $\{n_{k_l}\}$, such that 
\[
	v^{i,n_{k_l}}(x) \to v^i(x), \quad i=1,2,
\]
for all $x \in B$, and some $v^i:B \to \bR$. Given (\ref{grad}), we know that the sequence $\{v^{i,n}\}_{n \in \bN}$ is locally Lipschitz uniformly in $n$, and thus we can extend $\{v^i\}_{i=1,2}$ to the whole of $\bR^N$ by continuity. \textit{NB: For the sequel we use $n$ instead of $n_{k_l}$ for notational simplicity, assuming we work with this subsequence.}

We now prove that the sequence $\{ Z^{i,n} \}_{n \in \bN}$ is Cauchy. It is known (see, e.g. \cite{Hu}) that, since $Z^{i,n}$ is Markovian, it has the following representation:
\[
	Z^{i,n}_t = \nabla v^{i,n}(X_t) \sigma, \quad \text{ for all } t \geq 0.
\]
Hence, by (\ref{grad}), we see that 
\begin{equation}
	\bE \| Z^{i,n}_t \| \leq \bar{C} \bE (1 + \| X_t \|^2) \leq \tilde{C} (1 + \| x_0 \|^2),
\label{UZ}
\end{equation}
for some constants $\bar{C},\tilde{C}$ that are uniform in $t,n$. Now, set $Y^{i,n}_t = v^{i,n}(X_t)$ for $i=1,2$ and $n \in \bN$. Denote $\bar{Y} = Y^{i,n} - Y^{i,m}$, $\bar{Z} = Z^{i,n} - Z^{i,m}$, $\bar{\lambda}^i = \lambda^{i,n}-\lambda^{i,m}$ for $i=1,2$. For all $T \geq 0$, applying It\^o's lemma to $\bar{Y}^2$ and taking expectations, we obtain
\begin{equation}
	\bar{Y}_0^2 = \bE \big[ \bar{Y}_T^2 \big] - 2 \bE \bigg[ \int_0^T \bar{Y}_t \bar{f}(t)  dt \bigg] + \bar{\lambda}\bE \bigg[ \int_0^T \bar{Y}_t dt \bigg]+ \bE \bigg[ \int_0^T \| \bar{Z}_t \|^2 dt \bigg] ,
\label{EE2}
\end{equation}
where $\bar{f}(t) = f^{i,n}(X_t,Z^{i,n}) - f^{i,m}(X_t,Z^{i,m})$ in the notation above. Using the definition of $\bar{f}(t)$, estimate (\ref{grad}) and Remark \ref{remark2}, we see that
\[
	\bE | \bar{f}(t) |^2 \leq \bar{c} \bE (1 + \| X_t \|^4) \leq c(1 + \| x_0 \|^4),
\]
where $c$ and $\bar{c}$ are some constants that do not depend on $n,m$. Thus, applying Cauchy--Schwartz and the dominated convergence theorem, we conclude that 
\[
	\bE \bigg[ \int_0^T \bar{Y}_s \bar{f}(t)  dt \bigg] \to 0,
\]
as $n,m \to \infty$. The other terms in (\ref{EE2}) tend to zero by construction, and thus $\{ Z^{i,n} \}_{n \geq 0}$ is convergent in $\L^2_T(W)$ for all $T \geq 0$, where
\[
	\L^2_T(W) := \bigg \{  \text{predictable processes } \theta: \bE \bigg[ \int_0^T \| \theta_t \|^2dt \bigg] < \infty  \bigg \}.
\] 
In other words, there exist limits $Z^1, Z^2$, such that 
\[
	 \bE \bigg[ \int_0^T \| Z^{i,n}_t - Z^i_t\|^2 dt \bigg] \to 0, \quad i=1,2
\]
for all $T \geq 0$. 

It remains to show the driver term is also convergent in the appropriate topology. For all $T \geq 0$ we have 
\begin{equation}
\begin{split}
	&\bE \int_0^T \big|  f^{i,n}(X_s,Z^{i,n}_s) -   H_i(X_{s}, Z^{i}_s,u^*(Z^{1}_s,Z^{2}_s),v^*(Z^{1}_s,Z^{2}_s))  \big| ds \\
	       & \leq \bE \int_0^T  \big|  f^{i,n}(X_s,Z^{i,n}_s) -   f^{i,n}(X_s,Z^{i,n-1}_s)  \big| ds \\
	       & \qquad + \bE \int_0^T \big|  f^{i,n}(X_s,Z^{i,n-1}_s) -   H_i(X_{s}, Z^{i}_s,u^*(Z^{1}_s,Z^{2}_s),v^*(Z^{1}_s,Z^{2}_s))  \big| ds.
\end{split}	
\label{E3}
\end{equation}
By definition of the Hamiltonian, 
\[
	 \big|  f^{i,n}(X_s,Z^{i,n}_s) -   f^{i,n}(X_s,Z^{i,n-1}_s)  \big| \leq C \big\| Z^{i,n}_s - Z^{i,n-1}_s  \big\|,
\]
showing the convergence of the first term in (\ref{E3}) to zero. We recall that convergence in $L^2$ implies convergence in measure, which in turn implies convergence $a.e.$ for a subsequence. Therefore there exists a subsequence $\{ n_k \}_{k \in \bN}$, such that 
\[
	Z^{i,n_k} \to Z^i \quad \bP\otimes dt-a.e., \quad i=1,2.
\]
Recall that 
\[
	f^{i,n}(X_s,Z^{i,n-1}_s) = H_i(X_{s}, Z^{i,n-1}_s,u^*(Z^{1,n-1}_s,Z^{2,n-1}_s),v^*(Z^{1,n-1}_s,Z^{2,n-1}_s))
\]
and, given the continuity assumption on $f^{i,n}(X_s,Z^{i,n-1}_s)$ in $(Z^{1,n-1}_s,Z^{2,n-1}_s)$, we conclude that 
\[
	f^{i,n_k}(X_s,Z^{i,n_{k-1}}_s) \to H_i(X_{s}, Z^{i}_s,u^*(Z^{1}_s,Z^{2}_s),v^*(Z^{1}_s,Z^{2}_s)) \quad \bP \otimes dt-a.e.
\]
and convergence of the second term in (\ref{E3}) follows by the dominated convergence theorem. We have shown that 
\[
	\lim_{n \to \infty} f^{i,n_k}(X_s,Z^{i,n_k}_s) = H_i(X_{t}, Z^{i}_s,u^*(Z^{1}_s,Z^{2}_s),v^*(Z^{1}_s,Z^{2}_s)) \quad \text{in } \L^2_T(\bR^N) 
\]
for all $T\geq 0$. This concludes the proof.

\end{proof}

\begin{remark} We can now see why part (ii) in Assumption \ref{Isaacs} was necessary. In order to use the BSDE approach, we had to prove that a solution to a certain system of equations exists. We constructed such solution as a limit, and then had to show that the driver terms converged as well. Therefore, continuity of the Hamiltonian (which played the role of the driver) in $(z_1,z_2)$ was naturally required. 

\end{remark}

\begin{remark} One interesting question is the connection between solutions to finite horizon games (i.e. systems of BSDEs with continuous coefficients) and their ergodic counterpart. For a detailed account of the one dimensional Lipschitz case, see Hu, Madec and Richou \cite{HMR}. In the present context, given that at the optimum the drivers are only continuous, a measure change technique is not available, so the techniques used in \cite{HMR} cannot be easily transferred. 

One simple line of attack would be to fix one player's strategy at the ergodic optimum and focus on the long-run behaviour of the other player. This reduces the problem to the known case,  but raises more questions than it answers, since the control is no longer closed loop. Thus, even though convergence of the solution is guaranteed (as in \cite{HMR}), the interpretation of this system of equations as a stochastic game is lost. 

Considering the system of BSDEs with continuous coefficients directly is a challenging task. In particular, as there are no uniqueness results in this context, and measure change techniques are not available, it proves difficult to demonstrate convergence of solutions without significant further exploration of this theory. We leave these questions for future research.
\end{remark}

\subsection{Quick detour into one dimensional EBSDEs} In this section we leave the world of games and control and apply the machinery established above to the case of one dimensional ergodic BSDE with a continuous driver of linear growth (in the $z$ component). We begin by proving an auxiliary lemma that establishes a useful representation of this class of functions. 
\begin{lemma} Suppose we are given a function $f:\bR^d \to \bR$, that is continuous and of linear growth, that is, there exists a constant $\kappa > 0$, such that $|f(x)| \leq \kappa (1 + \| x \|)$ for all $x \in \bR^N$. Then there exist uniformly bounded functions $\phi(\cdot),\psi(\cdot)$, such that 
\[
	f(x) = \phi(x)x + \psi(x),
\]
holds for all $x \in \bR^d$.  
\label{L2}
\end{lemma}
 \begin{proof} We first notice that, given that $f$ is of linear growth, there exists a constant $C>0$ such that $\| f(x) \| \leq C\|x\|$ for all $x:\| x \| \geq 1$. Define 
\[
	\phi(x) := \bone_{\| x \| \geq 1}\frac{f(x)}{\| x \|^2}x^\top, \quad \psi(x) := \bone_{\| x \| < 1}f(x),
\]
where $x^\top$ is the transpose of $x$. The claim follows immediately. 

\end{proof}

\begin{remark} Notice that we have no continuity assumption on the individual components $\phi,\psi$. One can clearly see that this decomposition is inspired by the structure of the Hamiltonian in the games framework. 
\end{remark}

 Taking into account the result of Lemma \ref{L2}, we are therefore interested in the equation 
\begin{equation}
	 Y_t = Y_T+\int_t^T[f(X_s,Z_s)-\lambda] ds - \int_t^T Z_s dW_s,	
\label{E7}
\end{equation}
where one can find functions $\phi(\cdot),\psi(\cdot)$, such that $f(x,z) = \phi(x,z)z + \psi(x,z)$, and there exist constants $c,C > 0$, such that 
\begin{equation}
	\| \phi(x,z) \| \leq c, \quad | \psi(x,z) | \leq C,
\label{E6}
\end{equation}
holds for all $x,z \in \bR^N$. In general, the existence of a solution to (\ref{E6}) does not follow from the standard theorems about EBSDEs, for the driver $f(\cdot,\cdot)$ may not be Lipschitz. However, the technique developed for the game setup still gives the following result: 

\begin{theorem} Suppose the driver $f:\bR^N\times \bR^N \to \bR$ is a measurable map, such that $f(x,\cdot)$ is continuous for all $x \in \bR^N$ and there exists a constant $\kappa > 0$, such that 
\[
	 \big| f(x,z) \big| \leq \kappa (1 + \| z \|).
\] 
Then there exists a solution $(Y,Z,\lambda)$ to the EBSDE (\ref{E7}). Moreover, there exists a deterministic locally Lipschitz function $v:\bR^N \to \bR$, such that $Y_t = v(X_t)$ for all $t \geq 0$. 
\end{theorem}

\begin{proof} We use the iteration method developed in the previous section in the games context. Set $Z^0_s \equiv 0$, and let $(Y^n,Z^n,\lambda^n)$ be a solution to 
\begin{equation}	
	 Y^n_t = Y^n_T+\int_t^T[ \phi(X_s,Z^{n-1}_s)Z^n_s + \psi(X_s,Z^{n-1}_s)-\lambda^n] ds - \int_t^T Z^{n}_s dW_s,
\label{E7a}
\end{equation}
where the process $\{ Z^{n-1}_s \}_{s \geq 0}$ is a known Markovian solution we obtained on the previous iteration. In other words, $Z^{n-1}_s = \zeta^{n-1}(X_s)$ for some deterministic function $\zeta^{n-1}:\bR^N \to \bR^N$. In the notation of Lemma \ref{L2}, set 
\[
	f^n(x,z) = \phi(x,\zeta^{n-1}(x))z + \psi(x,\zeta^{n-1}(x)).
\]
We notice that there exist constants $c,C>0$, such that $| f^n(x,0) | \leq C$, and $| f^n(x,z) - f^n(x,z') | \leq c \| z - z' \|$, so by Theorem \ref{TErg} there exists a solution $(Y^n,Z^n,\lambda^n)$ to (\ref{E7a}). There also exist deterministic functions $y^n:\bR^N \to \bR, \zeta^n:\bR^N \to \bR^N$, such that $Y^n_s = v^n(X_s), Z^n_s = \zeta^n(X_s)$. Moreover, estimates (\ref{est}) and (\ref{grad}) hold uniformly in $n$. We can therefore use the diagonalisation  procedure on some dense countable subset $B \subset \bR^N$, in order to obtain the function $v:B \to \bR$, and $\lambda \in \bR$, such that 
\[
	v(x) = \lim_{k \to \infty }v^{n_k}(x) \quad \text{for all } x \in B, \quad \text{and } \lambda = \lim_{k \to \infty} \lambda^{n_k} 
\]
for some subsequence $\{ n_k \}_{k \geq 0}$. We can then extend $v$ to the entirety of $\bR^N$ by continuity. We then proceed to show convergence of $\{ Z^{n_k} \}_{k\geq 0}$ in $\L^2_T(W)$ for all $T \geq 0$ in exactly the same way as in the proof of Theorem \ref{T4}. We finally extract a further subsequence that converges $\bP \otimes dt - a.e.$, and use the fact that 
\[
\begin{split}
	| f^n(X_s,Z^n_s) &- f(X_s,Z_s) | 
		\\ & \leq | f^n(X_s,Z^n_s) - f^n(X_s,Z^{n-1}_s) | + | f^n(X_s,Z^{n-1}_s) - f(X_s,Z_s) |
		\\ & \leq c \| Z^n_s - Z^{n-1}_s \| + | f(X_s,Z^{n-1}_s) - f(X_s,Z_s) |.
\end{split}
\]
Using the dominated convergence theorem along with the continuity of $f$, we conclude the convergence of the driver term. 

\end{proof}

\begin{remark} We notice that the proof above extends naturally to the case of multiple dimensional EBSDEs with an appropriate diagonal boundedness structure. In other words, consider the system 
\[
	 Y^i_t = Y^i_T+\int_t^T[f^i(X_s,\bZ_s)-\lambda] ds - \int_t^T Z^i_s dW_s,	\quad i=1,\dots,n,
\]
where $\bZ$ denotes the vector of $Z$-components. Suppose for all $i = 1,\dots,n$ and for all $x \in \bR^N$, the map $\bz \to f^i(x,\bz)$ is continuous, and for every $i$, we know $ |f^i(x,\bz)|\leq \kappa(1+\|z^i\|)$ for some $\kappa>0$. Then the system admits a solution. 
\end{remark}

\begin{remark} Similar to the setting of finite horizon BSDEs with continuous coefficients considered by Lepeltier and San Martin (see \cite{LM}), we do not obtain any uniqueness result. Of course one could only hope for the uniqueness of the ergodic values, since $Y+c$ for any constant $c$ also solves the same equation. The difficulty in applying the standard technique to establish uniqueness for $\lambda$'s is that this uses a change of measure to remove a drift term of the form $[f(X_s,Z_s)-f(X_s,Z'_s)]ds$. This is not possible in the present framework due to the fact that $f(x,\cdot)$ is not Lipschitz. 
\end{remark}
%%%%%%%%%%%%UNIQUENESS%%%%%%%%%%%%%%%%%%

\subsection{Connection to PDEs}

In this section we briefly show the relation between Theorem \ref{T4} and the existence of a solution to the following system of elliptic ergodic HJB equations:
\begin{equation}
\begin{cases}
	\L y^1(x) + H_1(x,\nabla y^1(x) \sigma,u^*(x),v^*(x)) = \lambda^1, \\
	\L y^2(x) + H_2(x,\nabla y^2(x) \sigma,u^*(x),v^*(x)) = \lambda^2,
\end{cases}
\label{S2}
\end{equation}
where 
\[
	u^*(x) = u^*(\nabla y^1(x)\sigma,\nabla y^2(x)\sigma), \quad v^*(x) = v^*(\nabla y^1(x)\sigma,\nabla y^2(x)\sigma),
\]
the functions $y^i(\cdot)$, $i=1,2$ are as in Theorem \ref{T4}, and $\L$ is a linear operator defined as  
\[
	\L v(x) = \frac{1}{2}Tr \bigg( \sigma \sigma'(x) \nabla^2 v(x) \bigg) + \langle Ax + F(x) , \nabla v(x) \rangle.
\]
We notice that if we define the transition semigroup $\{ P_t \}_{t \geq 0}$ associated with the process $\{X_t\}_{t \geq 0}$ by 
\[
	P_t [\phi] (x) = \bE [\phi(X_t)],
\]
for all measurable functions $\phi:\bR^N \to \bR$ of polynomial growth, then $\L$ is the generator of $\{P_t\}_{t \geq 0}$. Similarly to Debussche, Hu and Tessitore \cite{Hu}, we adopt the following definition:
\begin{definition} Two pairs $\{ (y^i,\lambda^i) \}_{i=1,2}$ ($y^i:\bR^N \to \bR$, $\lambda^i \in \bR$) are called a mild solution to (\ref{S2}), if 
\begin{enumerate}[(i)] \item The functions $\{y^i(x)\}_{i=1,2}$ are in $C^1$ and of polynomial growth.
\item The gradient functions $\{ \nabla y^i (x)\}_{i=1,2}$ are of polynomial growth.
\item For all $0 \leq t \leq T$, for all $x \in \bR^N$, and for $i=1,2$,
\[
	y^i(x) = P_{T-t} [y^i](x) + \int_t^T \big( P_{T-s}\big[H_i(x,\nabla y^i(\cdot) \sigma,u^*(\cdot),v^*(\cdot))\big] (x) - \lambda^i \big)ds, 
\]
where $u^*,v^*$ are as above. 
\end{enumerate}
\end{definition}
 We therefore automatically have the following result: 
\begin{theorem} Suppose conditions in Theorem \ref{T4} are satisfied. Then $\{(y^i,\lambda^i)\}_{i=1,2}$ is a unique mild solution to the system (\ref{S2}). Conversely, if $\{(\y^i,\tilde{\lambda^i})\}_{i=1,2}$ is a mild solution to $(\ref{S2})$, then, setting $Y^i_t = \y^i(X_t)$, $i=1,2$, we obtain a solution to the system of EBSDEs (\ref{EBSDE}).
\end{theorem}

\color{black}

\subsection{Games with $n > 2$ players}

In this section we consider a game between $n \in \bN \cup \{\infty\}$ ergodic players. As in the case of two players, we show that the existence of a solution to a certain system of ergodic BSDEs guarantees the existence of a saddle point (provided a version of generalised Isaac's condition is assumed). We then proceed to prove the former. We show that the scenario of finitely many players is a relatively straightforward extension of Theorem \ref{T4}, but extra care is required for the countable case. 

 Let  $\bu = (u_1,\dots,u_n)$ be the vector of controls, such that $(\omega,t)\mapsto \bu_t(\omega)$ is a predictable process with values in a separable metric space $\U = \U_1 \times \dots \times \U_n$. Denote $\bz = (z_1,\dots,z_n) \in \bR^{n \times N}$. We define the Hamiltonian functions $\{H_i\}_{ i = 1,\dots,n}$, as 
\[
	H_i(x,z_i,\bu) = z_i R(\bu) + L_i(x,\bu),
\]
where, similar to the two player case, $L_i:\bR^N\times \U \to \bR$, $i=1\dots n$ are measurable cost functions Lipschitz in $x$, and the function $R:\U \to \bR^N$ corresponds to the control of the drift by players through the change of measure $\bP^{u,v,T} := \rho^{u,v}_T\bP$, where
\[
	\rho^{\bu}_T := \exp \bigg( \int_0^T R(\bu_t)dW_t - \frac{1}{2}\int_0^T \| R(\bu_t) \|^2 ds \bigg).	
\]
The payoffs are defined as 
\[
	J^i( \bu) = {\lim\sup}_{T\to\infty} T^{-1} \bE^{\bu,T}\bigg[\int_0^T L_i(X_t, \bu_t) dt\bigg],
\]
for $i=1\dots n$, and the goal is to find admissible vector of controls $\bu^*$, such that for all $i=1 \dots n$, 
\[
	J^i(u^*_i,\bu_{-i}) \geq J^i(\bu^*), 
\]
holds for all admissible controls $\bu_{-i}$, where $\bu_{-i}$ denotes the vector $\bu$ without its $i$-th component. The generalised Isaac's conditions become:
\begin{assumption} Suppose $n < \infty$. \begin{enumerate}[(i)]
 \item There exists a measurable map $\bu^*:\bR^{(n+1)N} \to \U$, such that for each $i = 1 \dots n$
\[
	H_i^*(x,\bz) := H_i(x,z_i,\bu^*(x,\bz)) \leq H_i(x,z_i,u_i,\bu_{-i}^*(x,\bz))
\]
 holds for all $(x,\bz,u_i) \in \bR^{(n+1)N} \times \U_i$. 
 \item The mapping $\bz \in \bR^{n\times N} \to \bH^*(x,\bz)$ is continuous for any fixed $x$. Here $\bH^*$ denotes the vector of $H^*_i$.
\end{enumerate}
\label{I2}
\end{assumption}

 As an immediate extension of the two player case we get the following result, the proof of which is identical to that of Theorem \ref{T3} with certain amount of notational care:
\begin{theorem} Given Assumption \ref{I2}, suppose that there exist $n \geq 2$ triplets $(Y^i, Z^i, \lambda^i)$, $i=1 \dots n$, such that 
\begin{equation}
	 Y^{i}_t = Y^i_T+\int_t^T[H_i(X_{s}, Z^i_s,\bu^*(\bZ_s))-\lambda^{i}] ds - \int_t^T Z^i_s dW_s,	
\label{ES}
\end{equation}
holds for all $0 \leq t \leq T < 0$. Moreover, assume there exist $n$ deterministic functions with polynomial growth $y_i(x)$, $i=1 \dots n$, such that $Y^i_t = y_i(X_t)$ holds $\bP$ - a.s. for all $t \geq 0$. Then the control $\bu^*(X_s,\bZ_s)$ is a Nash equilibrium.
\label{T5}
\end{theorem}
 For the case of finitely many players, we can show that the system (\ref{ES}) admits a solution using exactly the same Picard iteration as in the proof of Theorem \ref{T4}. In other words, we set $\bZ^0 \equiv 0$, and given a set of  Markovian solutions $\{(Y^{i,n-1},Z^{i,n-1},\lambda^i)\}_{i=1,\dots,n}$ obtained previously, we define the next iteration as a solution to 
\[
	 dY^{i,n}_t = - [f^{i,n}(X_t,Z^{i,n}_t)-\lambda^{i,n}]dt - \int_t^T Z^{i,n}_t dW_t,		\quad i=1,\dots,n,
\]
where 
\[
	f^{i,n}(X_s,Z^{i,n}) = H_i(X_{s}, Z^{i,n}_s,\bu^*(\bZ^{n-1}_s)). 
\]
The rest of the proof is identical.

 \subsubsection{Infinite player games} For a countable number of players additional care is needed. First, condition (ii) in Assumption \ref{I2} needs to be reformulated. This is due to the fact that in infinite dimensions not all norms are equivalent, and therefore in order to have a notion of continuity, we need to specify the topology on infinite dimensional vectors $\bz$. Let $\theta^i$ denote the $i$-th component of the vector $\theta$, and define
\[
	\L^2_{\epsilon,loc} := \bigg \{  \text{predictable processes } \theta: \bE \bigg[ \sum_{i=1}^{\infty} \frac{1}{i^{1+\epsilon}}\int_0^T\| \theta^i_t \|^2 dt \bigg] < \infty,  \text{ for all }T \bigg\}.
\]
Given the fact that constants in the estimate (\ref{UZ}) do not depend on time or the number of iteration, we know that there exists a constant $C > 0$ such that
\[
	\bE \bigg[ \int_0^T \| Z^{i,n} _t\| dt\bigg] < CT, \quad \text{ for all } T > 0,
\]
holds uniformly in $i,n$. Hence $\bZ^{n} \in \L^2_{\epsilon,loc}$ for all $n \geq 0$. It is clear that componentwise convergence of $\bZ^n$ to $\bZ$ in $\L^2_T(W)$ for all $T > 0$ guarantees convergence of $\bZ^n$ to $\bZ$ in $\L^2_{\epsilon, loc}$. We also notice that, in the proof of Theorem \ref{T4} we use continuity of the vector of Hamiltonians $\bH^*$ for each component separately. Therefore, the infinite dimensional analogue of (i) in Assumption \ref{I2} remains exactly the same, while (ii) becomes 
\begin{enumerate}[\textbf{(ii)}] \item The mapping $\bz \in \L^2_{\epsilon,loc} \to H_i^*(x,\bz)$ is continuous (for $\bP \times dt$ almost all $(\omega,t)$) for each $i \in \bN$. 
\end{enumerate}
Secondly, the construction of a convergent subsequence for the vector of ergodic values $\bl^{k}$ needs modification. This is due to the fact that infinite dimensional cubes are not compact. However, the result of Lemma \ref{L1} clearly holds, ensuring that each component $\lambda^{i,k}$ lies in a compact (uniformly in $k$), and thus we can use a diagonalisation procedure to construct a limit $\bl$. Using \textbf{(ii)} instead of (ii), the rest of the proof remains the same.

\begin{remark} Even though the topology $\L^2_{\epsilon,loc}$ looks artificial, it has an intuitive interpretation. As pointed out in Remark \ref{remarkTmp}, the generalised Isaac's condition is an equivalent of Nash equilibrium on the infinitesimal scale. Hence the Hamiltonians $\{H_i\}_{i \geq 1}$ in some sense represent ``local values'' for each player. We also note that the vector $\{\bZ_t\}_{t \geq 0}$ can be seen as a vector of strategies, since it determines the values of optimal controls $\{\bu^*_t\}_{t \geq 0}$. 

 Therefore continuity of $\{H_i\}_{i \geq 1}$ in $\bz \in \L^2_{\epsilon,loc}$ means that, when finitely many ``neighbours'' (the term can be made rigorous if we put a graph structure on the space of agents) of player $i \in \bN \cup \{ \infty \}$ perturb their strategies only slightly, the change in his local value will be insignificant as well, even if distant agents make substantial changes. 
\end{remark}

 Further generalisations may include the case of infinite dimensional forward process, the presence of L\'evy noise, as well as the introduction of time periodicity. Given that the corresponding existence and uniqueness results were established in \cite{Levy_paper} (for the case where $X$ takes values in a separable Hilbert spaces), we expect that all the results of the present paper should still hold true. 

 Given the existence result for a game between countably many players, the next natural step would be to look at ergodic mean-field type games. The difficulty with doing this is the continuity assumption on the Hamiltonian, which depends on the topology on $z$. Since our result is quite generic, it allows for the introduction of additional dependancies between policies of various players, provided that we still have the corresponding version of the generalised Isaac's condition.

\section{Stochastic games with asymmetric agents}
In this section we consider a game between two players of different type. The main goal is to understand whether there is a fundamental difference between an ergodic player (as defined above), and one who discounts the future as the discount rate becomes negligible. We define the payoff structure as follows:  
\[
\begin{split}
	&J^1(x_0, u,v) = {\lim\sup}_{T\to\infty} T^{-1} \bE^{u,v,T}\bigg[\int_0^T L_1(X_t, u_t,v_t) dt\bigg], \\
	&J^2(x_0, u,v) =  \bE^{u,v,T}\bigg[\int_0^{\infty} e^{-\alpha t}L_2(X_t, u_t,v_t) dt\bigg],
\end{split}
\]
where $\alpha$ is some positive constant. In the sequel we will be calling the player of the first type `ergodic', and of the second `discounting'. We begin by proving an analogue of Theorem \ref{T3} for this case:
\begin{theorem} Let the Hamiltonians $\{H_i\}_{i=1,2}$ be defined as in (\ref{H}) and suppose Assumption \ref{Isaacs} holds. Suppose further that there exist a triplet $(Y, Z, \lambda)$, and a pair $(\Y,\Z)$ such that 
\begin{equation}
\begin{cases}
	  Y_t = Y_T+\int_t^T[H_1(X_{s}, Z_s,u^*(Z_s,\Z_s),v^*(Z_s,\Z_s))-\lambda] ds - \int_t^T Z_s dW_s,	\\
	  \Y_t = \Y_T+\int_t^T[H_2(X_{s}, \Z_s,u^*(Z_s,\Z_s),v^*(Z_s,\Z_s))-\alpha \Y_s] ds - \int_t^T \Z_s dW_s,
\end{cases}
\label{S1}
\end{equation}
holds for all $0 \leq t \leq T < 0$. Moreover, assume there exist two deterministic functions with polynomial growth $y(x),\y(x)$, such that $Y_t = y(X_t)$ and $\Y_t=\y(X_t)$ holds $\bP$-a.s. for all $t \geq 0$. Then the control $(u^*(X_s,Z_s,\Z_s),v^*(X_s,Z_s,\Z_s))$ is a Nash equilibrium.
\label{T6}
\end{theorem}

 \begin{proof} We first notice that, since the Hamiltonian structure remains the same, the proof for the ergodic player is identical to that in Theorem \ref{T3}. We therefore focus on the discounting one. For an arbitrary admissible control $v$, the pair $(u^*,v)$ generates a measure $\bP^{u^*,v}$ under which $dW^{u^*,v}_t = dW - R(u^*_t,v_t)dt$ is a Brownian motion. Then 
\[
	d (e^{-\alpha t} \Y_t) = -e^{-\alpha t} \big( [ H_2(X_{t}, \Z_t,u^*_t,v^*_t) + \Z_t R(u^*_t,v_t) ] dt +  \Z_t dW^{u,v^*}_t \big),
\]
and thus, integrating and taking expectations with respect to $\bP^{u^*,v}$, we obtain
\[
\begin{split}
	\Y_0 &=  \bE^{u^*,v} \bigg[ \int_0^T \big(  H_2(X_{t}, \Z_t,u^*_t,v^*_t) - \Z_t R(u^*_t,v_t) - L_2(X_t,u^*_t,v_t)  \big)dt \bigg]
		\\& \qquad + e^{-\alpha T}\bE^{u^*,v} \big[ \Y_T \big]+  \bE^{u^*,v} \bigg[ \int_0^T L_2(X_t,u^*_t,v_t) dt\bigg].
\end{split}
\]
Taking the $\limsup$ on both sides, we obtain 
\[
	\Y_0 \leq \bE^{u^*,v} \bigg[ \int_0^{\infty} L_2(X_t,u^*_t,v_t) dt\bigg] = J^2(x_0,u^*,v),
\]
and the equality holds for $v = v^*$. 

\end{proof}

\begin{theorem} The system (\ref{S1}) admits a Markovian solution in the sense of Theorem \ref{T6}.
\end{theorem}
 \textbf{Sketch of the proof:} Since the rigorous proof is almost identical to that of Theorem \ref{T4}, we only outline the minor differences. Our first aim is to construct a solution to the system (\ref{S1}) using Picard iteration (\ref{E2}). Notice that for a fixed $\alpha > 0$, the uniform estimates (\ref{est}) and (\ref{grad}) hold for the ergodic BSDE as before, while for the discounted one, (\ref{est}) is replaced by 
\[
	| v^{2,n}(x) | \leq C/\alpha, \quad | v^{2,n}(x) - v^{2,n}(0)| \leq C(1 + \| x \|^2),
\]
where the constant $C>0$ is independent of $n$ and $\alpha$ (for details see Theorem 8 in \cite{Levy_paper}). This allows us to follow the rest of the proof of Theorem \ref{T4} with only slight modifications. 

\qed

 We are now interested in what happens when we send the discount rate $\alpha$ to zero. The idea is to extract a convergent subsequence from a sequence $\{ \alpha_n \}_{n \in \bN}$ and then study the limiting object. We begin by showing the following result:

\begin{theorem} Let $y^{\alpha}(x),\y^{\alpha}(x)$ be the $Y$ components of the Markovian solution to (\ref{S1}), and $\lambda^{\alpha}$ be the ergodic value of the first player. Then there exists a sequence $\alpha_n \downarrow 0$ and two functions of polynomial growth $y(\cdot)$ and $\y(\cdot)$, such that 
\[
	y^{\alpha_n}(x) \to y(x), \quad \lambda^{\alpha_n} \to \lambda, \quad \big(\y^{\alpha_n}(x) - \y^{\alpha_n}(0)\big) \to \y(x), \quad \alpha_n \y^{\alpha_n}(0) \to \tilde{\lambda},
\]
holds for all $x \in \bR^N$.
\end{theorem}
 \begin{proof} We know that, by construction, the functions $y^{\alpha_n}(\cdot)$ and $\y^{\alpha_n}(\cdot)$ are locally Lipschitz, uniformly in $n$. Also,
\[
	| y^{\alpha_n}(x) | \leq C(1 + \| x \|^2), \quad | \y^{\alpha_n}(x) - \y^{\alpha_n}(0) | \leq C(1 + \|x \|^2)
\]
holds for all $x \in \bR^N$. By Lemma \ref{L1} we also know that $\{\lambda^{\alpha_n}\}_{n \geq 1}$ lie in a compact. We can therefore use the diagonalisation  procedure to construct the limit on a dense subset of $\bR^N$ (e.g. $\bQ^N$), and then extend by continuity to the entire space. 

\end{proof}

 Using arguments identical to the those in the proof of Theorem \ref{T4}, we conclude that as the discount rate of the second player goes to zero, we can extract a subsequence such that 
\[
	J^{1,\alpha_n}(x_0,u^*,v^*) = \lambda^{\alpha_n} \to \lambda, \quad \alpha_n J^{2,\alpha}(x_0, u^*,v^*) = \alpha_n \y^{\alpha_n}(x_0) \to \tilde{\lambda}, 
\]
and we can find $(\lambda,\tilde{\lambda})$ by solving the following system:
\[\begin{cases}
	  Y_t = Y_T+\int_t^T[H_1(X_{s}, Z_s,u^*(Z_s,\Z_s),v^*(Z_s,\Z_s))-\lambda] ds - \int_t^T Z_s dW_s,	\\
	  \Y_t = \Y_T+\int_t^T[H_2(X_{s}, \Z_s,u^*(Z_s,\Z_s),v^*(Z_s,\Z_s))-\tilde{\lambda}] ds - \int_t^T \Z_s dW_s.
\end{cases}
\]
This is in line with our intuition. Indeed, we know that for the case of one player (i.e. a standard optimal control problem), the EBSDE is obtained a a limit of discounted infinite horizon BSDEs.

\bibliographystyle{plain}
\bibliography{Bibliography}

\end{document}